# QUANTUM ERGODICITY OF BOUNDARY VALUES OF EIGENFUNCTIONS: A CONTROL THEORY APPROACH

N. BURQ

ABSTRACT. Consider $M$, a bounded domain in $\mathbb{R}^d$, which is a Riemanian manifold with piecewise smooth boundary and suppose that the billiard associated to the geodesic flow reflecting on the boundary acording to the laws of geometric optics is ergodic. We prove that the boundary value of the eigenfunctions of the Laplace operator with reasonable boundary conditions are asymptotically equidistributed in the boundary, extending previous results by Gérard, Leichtnam [7] and Hassel, Zelditch [8] obtained under the additional assumption of the convexity of $M$.

RÉSUMÉ. Soit $M$ un domain borné de $\mathbb{R}^d$ qui est une variété riemanienne à coins. On suppose que le billard défini par le flot géodésique brisé est ergodique. On démontre que les valeurs au bord des fonctions propres du Laplacien (avec des conditions aux limites raisonnables) sont asymptotiquement équidistribuées dans le bord. Ceci généralise des résultats antérieurs de P. Gérard, E. Leichtnam [7] et A. Hassel, S. Zelditch [8] obtenus sous l'hypothèse supplémentaire de convexité géodésique du domaine.

## 1. INTRODUCTION

The purpose of this article is to show how ideas coming from control theory allow to deduce ergodicity of boundary values of eigenfunctions on manifolds with piecewise smooth boundaries and ergodic billiard flows from the known results on ergodicity of the eigenfunctions in the interior. The main result we obtain is the following:

**Theorem 1.** *Consider the eigenfunctions of the Laplace operator associated to the boundary condition B:*

(1.1)
$$-\Delta_g e_j = \lambda_j^2 e_j$$
$$Be_j \mid_{\partial M} = 0$$

*Denote by $|dx'|$ and $d\sigma$ the natural measures on $\partial M$ and $T^*\partial M$; and associate to $e_j$ its boundary value according to the table below*

| $B$ | $Be$ | $e^b$ | $d\mu_b$ |
|---|---|---|---|
| $Dirichlet$ | $e \mid_{\partial M}$ | $\lambda_k^{-1}\partial_n e \mid_{\partial M}$ | $\gamma(\varrho)1_{\|\varrho\|_{g(x)}\leq 1}d\sigma$ |
| $Neumann$ | $\partial_n e \mid_{\partial M}$ | $e \mid_{\partial M}$ | $\gamma^{-1}(\varrho)1_{\|\varrho\|_{g(x)}\leq 1}d\sigma$ |
| $Robin$ | $(\partial_n e - Ke) \mid_{\partial M}$ | $e \mid_{\partial M}$ | $\frac{\gamma(\varrho)}{\gamma^2(\varrho)+k(\varrho)^2}1_{\|\varrho\|_{g(x)}\leq 1}d\sigma$ |

1991 *Mathematics Subject Classification*. 35Q55, 35BXX, 37K05, 37L50, 81Q20 .

This work was carried out during a stay of the author at the University of California, Berkeley, funded by the Miller Institute for Basic Research in Science; and a stay at the Mathematical Science Research Institute, Berkeley. I thank these institutions for their support.

where $K \in \Psi^1(\partial M)$ is a non-negative self adjoint pseudodifferential operator of order 1 (this assumption could be relaxed to the non negativness of its principal symbol). Then there exists a subset $S$ of density $1$ such that for any semi-classical pseudodifferential operator on the boundary (tangential to the boundary of the boundary) of order $0$, $a(x, hD_x)$ we have

$$(1.2) \qquad \lim_{i \to +\infty; i \in S} \int_{\partial M} a(x, \lambda_j^{-1} D_x) e_j^b \overline{e_j^b} d\sigma = \langle \mu_b, a(x, \xi) \rangle.$$

The study of ergodicity of the eigenfunctions of the Laplace operator dates back to the works by Shnirelman [10], Zelditch [11] and Colin de Verdière [6]. In the case of manifold with boundary, the first study was initiated by Gérard and Leichtnam [7] for $C^{1,1}$ boundaries and then extended (with a much simpler proof) by Zelditch and Zworski [12] to the case of manifolds with corners.

Recently, the study of the ergodicity of the boundary values of the eigenfunctions has been performed in the restricted case of convex manifolds (with corners) by Hassel and Zelditch [8]. Theorem 1 is hence a generalization of these results in [8] to the non convex case. Since most known ergodic billiards are non convex, our generalization is relevant (see the discussion in the introduction of [8]).

In fact it has to be noticed that for Dirichlet boundary conditions and the standard metric, the result above is not new (even though it does not appear explicitly in the literature for non convex domains): indeed it had been recognized by Gérard and Leichtnam [7, Remark 5.3] that for the Dirichlet problem, the questions of ergodicity of the eigenfunctions in the interior and ergodicity of their normal derivatives on the boundary are equivalent. The assumption of convexity of the domain which was made throughout the article played no role in this particular remark; and to obtain the ergodicity of the normal derivatives for a general ergodic billiard, it is consequently sufficient to link the posterior result in [12] and this remark (and to treat some minor problems arising from the presence of the corners, see the end of Section 3).

On the other hand, in control theory, a lot of analysis has been done since the 80's to relate the behaviour of solutions of wave equations in the interior of a domain with the behaviour of their traces on the boundary (see [1, 2, 3]). We are going to take advantage of this knowledge (eigenfunctions are -very- particular cases of solutions of wave equations) in order to show that the interior result in [12] implies the boundary result. For general boundary conditions, the analysis becomes more involved than in the Dirichlet case; but the knowledge of [12] (the eigenfunctions are, for a density-1 sequence, equidistributed and consequently can not concentrate at the boundary) allows to avoid most of the technical difficulties (see [3]) usually posed by boundary value problems; and we shall be able to perform quite the same argument as in [7] in the case of boundary conditions satisfying the strong Lopatinskii condition (non degenerate Robin conditions for example). In the case of a Neumann problem, this strong Lopatinskii condition is no more fulfilled at the gliding points (and in particular the traces of the eigenfunctions are not *a priori* bounded in $L^2(\partial M)$) and one has to show in addition that nothing happens at these points. In fact this part of the analysis, as well as the study of the behaviour near the singular part of the boundary, will be performed using some heat kernel methods from [8, Appendix].

In section 2, we give the assumptions and define the semiclassical measures associated to sequences of eigenfunctions. We also give some elementary consequences of Zelditch and Zworski's result. We present an elementary and self contained (assuming the result in [12]) proof of Theorem 1 for Dirichlet conditions in section 3. In section 4 we prove the result for Robin and Neumann conditions for a restricted class of operators whose principal symbol vanishes near the singular set $\Sigma$ and the Glancing set $\mathcal{G}$. Finaly in section 5, we explain how the analysis in [8, Appendix] allows to handle these singular and glancing sets.

**Remark 1.1.** *After this work was completed, S. Zelditch informed me that they had simultaneously been able to supress the convexity assumption in their work with A. Hassel (see [8, Addendum]).*

## 2. Assumptions and ergodicity of eigenfunctions

The manifold $(M, g)$ we consider is a compact $C^\infty$ Riemann manifold with piecewise boundary: $M$ is a compact subset of a $C^\infty$ compact Riemann manifold, $\widetilde{M}$ and there exist $r$ functions $f_1, \ldots, f_r \in C^\infty(\widetilde{M})$ such that

- 1) $df_i \mid f_i^{-1}(0) \neq 0$,
- 2) $df_i$ and $df_j$ are independant on $f_i^{-1}(0) \cap f_j^{-1}(0)$,
- 3) $M = \{x \in \widetilde{M}; f_i(x) > 0 \forall i = 1, \ldots, r\}$

To include Bunimovitch's stadium example, condition 2) can be replaced by

- 2') $f_i^{-1}(0) \cap f_j^{-1}(0)$ is an embedded submanifold of $\widetilde{M}$ and $M$ has Lipschitz boundary.

A point $x_0 \in \partial M$ will be called regular if

$$\sharp\{j; f_j(x_0) = 0\} = 1 \tag{2.1}$$

and singular otherwise. We will denote by $\Sigma \subset \partial M$ the set of singular points.

### 2.1. Ergodicity of eigenfunctions.
With these assumptions, one can define the billiard flow (see [12]) and under the assumption that this flow is ergodic, S. Zelditch and M. Zworski prove the following

**Theorem 2** ([12]). *There exist a set $S \subset \mathbb{N}$ of density one:*

$$\lim_{N \to +\infty} \frac{1}{N} \sharp(S \cap \{1, 2, \ldots, N\}) = 1 \tag{2.2}$$

*such that for any $A \in \Psi^0_{phg}$ pseudodifferential operator of order $0$, with kernel supported in the interior of $M$, one has*

$$\lim_{j \in S \to +\infty} (Ae_j, e_j)_{L^2(M)} = \langle dL, \sigma_0(A) \rangle \tag{2.3}$$

*where $dL$ is the Liouville measure on $S^*M$ the sphere cotangent bundle.*

**Remark 2.1.** *In [12], this result is proved only for Dirichlet or Neumann boundary conditions. However, the proof relying on the easy (non-diffractive) part of Egorov theorem extends clearly to Robin boundary conditions.*

## 2.2. Semi-classical measures.

We introduce the following usual class of semi-classical symbols on $\mathbb{R}^d$:

(2.4) $$S^m(\mathbb{R}^d) = \{a \in C^\infty(T^*\mathbb{R}^d \times (0,1]); |\partial_x^\alpha \partial_\xi^\beta a(x,\xi,h)| \leq C_{\alpha,\beta,k} h^m \langle \xi \rangle^{-|\beta|}\}$$

To any symbol $a \in S^m$ we can associate an operator $\mathrm{Op}(a)(x, hD_x, h)$ by the relation

(2.5) $$\mathrm{Op}(a)(x, hD_x, h)u = \int e^{ix\cdot\xi} a(x, h\xi, h)\widehat{u}(\xi)d\xi.$$

Denote by $\widetilde{M}$ the space $\mathbb{R}^d$ endowed with a metric equal to $g$ on $M$. Consider the sequence of eigenfunctions $(e_k)_{k\in\mathbb{N}}$ normalized in $L^2(M)$; and denote by $\underline{e_k}$ their extension by 0 outside $M$ which is bounded in $L^2(\widetilde{M})$. Using the standard definition of semi-classical measures in [7, section 2], [9] (see also [4]) we can associate, modudo extraction of a subsequence, a measure $\mu$ on $T^*\widetilde{M}$ such that for any semiclassical pseudodifferential operator $\mathrm{Op}(a)(x, hD_x)$ one has (with $h_k = \lambda_k^{-1}$),

(2.6) $$\lim_{k\to+\infty} \left(a(x, h_k D_x)\underline{e_k}, \underline{e_k}\right)_{L^2(\widetilde{M})} = \langle \mu, a(x, \xi, 0)\rangle.$$

Of course this measure is supported in $T^*\widetilde{M}|_{\overline{M}}$. The semi-classical elliptic regularity implies that *in the interior of $M$*, the semi-classical measure $\mu$ is supported in the characteristic set:

(2.7) $$\mathrm{supp}(\mu) \subset \mathrm{Char}(-h^2\Delta_g - 1) = \{(x,\xi) \in T^*\mathbb{R}^d; |\xi|_{g(x)} = 1\}$$

The ergodicity property (2.3) and (2.7) imply that in the interior of $M$, $\mu$ satisfies

(2.8) $$\mu 1_{x\in M} = dL \otimes \delta_{|\xi|_{g(x)}=1}.$$

But due to the normalization $\|\underline{e_i}\|_{L^2} = 1$, we know that $\mu$ has at most total mass 1. Since it is also the case for $dL \otimes \delta_{|\xi|_{g(x)}=1}$, we conclude

(2.9) $$\mu = \mu 1_{x\in M} = dL \otimes 1_{|\xi|_{g(x)}=1}$$

and consequently for any semi-classical pseudodifferential operator on $\mathbb{R}^d$, one has

(2.10) $$\lim_{j\in S\to+\infty} \left(\mathrm{Op}(a)(x, h_j D_x)\underline{e_j}, \underline{e_j}\right)_{L^2(\widetilde{M})} = \langle dL \otimes \delta_{|\xi|_{g(x)}=1}, \mathrm{Op}(a)(x,\xi)\rangle.$$

Finally, we can extend slightly the class of symbols for which (2.5) holds. We work close to a point $x_0 \in \partial M \setminus \Sigma$ in a geodesic coordinate system $(x_n, x')$ such that locally $M = \{(x_n, x') \in \mathbb{R}^{*,+} \times \mathbb{R}^{d-1}\}$ and

(2.11) $$\Delta_g = \frac{1}{\sqrt{c(x)}}\partial_{x_n}\sqrt{c(x)}\partial_{x_n} + R(x_n, x', D_{x'})$$

and the natural integration measures on $M$ and $\partial M$ are

(2.12) $$|dx| = \sqrt{c(x)}dx_n dx', \qquad |dx'| = \sqrt{c(x)}dx'$$

In this coordinate system, the usual elliptic ($\mathcal{E}$), hyperbolic ($\mathcal{H}$) and glancing ($\mathcal{G}$) subsets of $T^*\partial M$ are

(2.13) $$\begin{aligned}\mathcal{E} &=, \{(x',\xi'); R(0,x',\xi' > 1\} = \{(x',\xi'); \|\xi'\|_{g(x)} > 1\} \\ \mathcal{G} &=, \{(x',\xi'); R(0,x',\xi' = 1\} = \{(x',\xi'); \|\xi'\|_{g(x)} = 1\} \\ \mathcal{H} &=, \{(x',\xi'); R(0,x',\xi' < 1\} = \{(x',\xi'); \|\xi'\|_{g(x)} < 1\}\end{aligned}$$

**Lemma 2.2.** *Denote by $S^m(\mathbb{R}_{x_n} \times T^*\mathbb{R}^{d-1}_{x'})$ the set of symbols in the $x', \xi'$ variable depending smoothly of the parameter $x_n$. For any $a_1, a_2 \in S^0(\mathbb{R}_{x_n} \times T^*\mathbb{R}^{d-1}_{x'})$*

$$\lim_{k \in S \to +\infty} (Op(a_1)(x_n, x', h_k D_{x'}, h_k) e_k, e_k)_{L^2(M)} = \langle \mu, a_1(x, \xi') \rangle. \tag{2.14}$$

$$\lim_{k \in S \to +\infty} \left(Op(a_2)(x_n, x', h_k D_{x'}) h_k D_{x_n} e_k, e_k\right)_{L^2(M)} = \langle \mu, a_2(x, \xi') \xi_n \rangle. \tag{2.15}$$

Indeed if the operators $a_1$ and $a_2$ have kernels supported in $M$, this result is straightforward (using the elliptic regularity). To allow symbols supported up to the boundary, remark first that if $\varphi \in C_0^\infty(\mathbb{R})$ is equal to 1 close to 0, due to (2.9) we have

$$\lim_{\varepsilon \to 0} \lim_{j \in S \to +\infty} \|\varphi(x_n/\varepsilon) e_j\|^2_{L^2(M)} = \lim_{\varepsilon \to 0} \langle \mu, \varphi(\cdot/\varepsilon)^2 \rangle = 0. \tag{2.16}$$

On the other hand

$$\begin{aligned}
0 &= \int_M (-h_j^2 \Delta_g - 1) e_j \varphi^2(x/\varepsilon) \overline{e_j} |dx| \\
&= \int_M \varphi^2(x/\varepsilon) h_j \nabla_g e_j \overline{h_j \nabla_g e_j} |dx| - \int_M \varphi^2(x/\varepsilon) e_j \overline{e_j} |dx| \\
&\quad + h_j \int_M \nabla_g(\varphi^2(x/\varepsilon)) h_j \nabla_g e_j \overline{e_j} |dx| + \int_{\partial M} h^2 K(e_j) \overline{e_j} |dx'|
\end{aligned} \tag{2.17}$$

(with $K = 0$ in case of Dirichlet or Neumann boundary conditions) which implies (using that the operator $K$ is non-negative)

$$\lim_{j \to +\infty} \int_M \varphi^2(x/\varepsilon) h_j \nabla_g e_j \overline{h_j \nabla_g e_j} |dx| \leq \lim_{j \to +\infty} \int_M \varphi^2(x/\varepsilon) e_j \overline{e_j} |dx| + \lim_{j \to +\infty} \mathcal{O}_\varepsilon(h_j) \tag{2.18}$$

and consequently

$$\lim_{\varepsilon \to 0} \lim_{j \to +\infty} \int_M \varphi^2(x/\varepsilon) h_j \nabla_g e_j \overline{h_j \nabla_g e_j} |dx| = 0 \tag{2.19}$$

Now (2.16), (2.19) (and the knowledge of the result if the operators are supported in the interior) imply Lemma 2.2.

## 3. The Dirichlet case

**3.1. A priori estimates.** We start this section by recalling some well known estimates on the boundary values of the normal derivatives of the eigenfunctions (see for example in the more general context of strong Lopatinskii conditions [5, Theorem VII.4.4] or in our context [7, Lemma 2.1]).

**Proposition 3.1.** *Consider $(e_j)$ the Dirichlet eigenfunctions of the Laplace operator associated to the eigenvalues $\lambda_j$ and normalized to 1 in $L^2(M)$. Denote by $e_j^b = \lambda_j^{-1} \partial_n e_j |_{\partial M}$ the normal derivatives on the boundary of the eigenfunctions. Then there exists $C > 0$ such that*

$$\|e_j^b\|_{L^2(\partial M)} \leq C \tag{3.1}$$

*Proof.* Consider $X$ a vector field defined near $\overline{M}$ and such that on $\partial M$ one has $X \cdot \partial_n \geq c > 0$, and compute

$$\int_M [-h_k^2 \Delta - 1, X] e_k \overline{e_k} dx = \int_M -h_k^2 \Delta X e_k \overline{e_k} |dx|,$$

(3.2)
$$= h_k^2 \int_{\partial M} X e_k \mid_{\partial M} \frac{\partial \overline{e_k}}{\partial n} |dx'|,$$

$$\geq c \int_{\partial M} \left| \frac{\partial e_k}{\partial n} \right|^2 |dx'|.$$

On the other hand, the left hand term in (3.1) has the form

(3.3)
$$\int_M h_k^2 P_2 e_k \overline{e_k}$$

where $P_2$ is a second order differential operator, and hence this quantity is clearly bounded by $h_k^2 \|e_k\|_{H_0^1(M)}^2 \leq C$.

### 3.2. The relationship between the interior and the boundary.
Consider a point $x_0 \in \partial M \setminus \Sigma$. Working in a geodesic coordinate system close to $x_0$, we can compute for symbols $a = a_1(x, \xi') + a_2(x, \xi') \xi_n$, using Green's formula

(3.4) $$\int_M \frac{1}{ih_k} [-h_k^2 \Delta_g - 1, \mathrm{Op}(a)(x, h_k D_x)] e_j \overline{e_j} |dx|$$

$$= \int_M \frac{1}{ih_k} (-h_k^2 \Delta_g - 1) \mathrm{Op}(a)(x, h_k D_x) e_j \overline{e_j} |dx'|$$

$$= \int_{\partial M} ih_k \partial_{x_n} \mathrm{Op}(a)(x, h_k D_x) e_j \overline{e_j} |dx'|$$

$$- \int_{\partial M} i \mathrm{Op}(a)(x, h_k D_x) \overline{e}_j \mid_{x_n=0} \overline{h_k \partial_{x_n} e_j} |dx'|$$

$$= \int_{\partial M} \mathrm{OpOp}(a_2)(x_n = 0, x', h_j D_{x'}) h_j \partial_n e_j \mid_{x_n = 0} \overline{h_j \partial_{x_n} e_j} |dx'|$$

Now consider the left hand side in (3.4)
(3.5)
$$\int_M \frac{1}{ih_k} [-h_k^2 \Delta_g - 1, \mathrm{Op}(a)(x, h_k D_x)] e_j \overline{e_j} |dx| = \int_M \mathrm{Op}(q)(x_n, x', h D_{x_n}, h D_{x'}) e_j \overline{e_j} |dx| + \mathcal{O}(h_j)$$

where

(3.6) $$q(x_n, x', \xi_n, \xi') = \{\xi_n^2 + R(x_n, x', \xi') - 1, a(x, \xi)\}$$

is a second order polynomial in the $\xi_n$ variable. Taking into account the equation $(-h^2 \Delta_g - 1)e = 0$ and (2.11), we can eliminate the second order term in $\xi_n$ and obtain (using (2.14), (2.15)) that the limit $j \in S \to +\infty$ of the left hand side in (3.4) is equal to

(3.7) $$\langle \mu, \widetilde{q} \rangle$$

where $\widetilde{q}$ is the first order polynomial in the $\xi_n$ variable equal to $q$ on the set $\{\xi_n^2+R(x_n,x',\xi')-1=0\}$. But the measure $\mu$ is supported in this set and consequently we get

$$\langle \mu, \widetilde{q}\rangle = \langle \mu, q\rangle \tag{3.8}$$

which implies that the left hand side in (3.4) tends to

$$-\langle \mu, H_p(a)\rangle = \langle H_p(\mu), a\rangle \tag{3.9}$$

where $p(x,\xi) = \xi_n^2 + R(x,\xi') - 1$. But since we know that the measure $\mu$ is equidistributed, we have far from the boundary $H_p(\mu) = 0$; and close to the boundary we can compute $H_p(\mu)$ by the jump formula:
(3.10)
$$H_p(\mu) = (2\xi_n \partial_{x_n}(1_{x_n>0})\mu = 2\sqrt{1-R}\delta_{x_n=0}\otimes\delta_{\xi_n=\sqrt{1-R}}\otimes|dx'| - 2\sqrt{1-R}\delta_{x_n=0}\otimes\delta_{\xi_n=-\sqrt{1-R}}\otimes|dx'|$$

and we obtain using (3.4) and (3.9)

$$\lim_{j\in S\to+\infty} \int_{\partial M} a_2(x_n=0,x',h_j D_{x'})h_j\partial_n e_j\,|_{x_n=0}\,\overline{h_j\partial_{x_n} e_j}|dx'| \tag{3.11}$$
$$= \int_{T^*\partial M} a_2\sqrt{1-R}1_{R\leq 1}|dx'||d\xi'| = \langle a_2\sqrt{1-R}, |dx|\otimes 1_{R\leq 1}\rangle$$

which implies Theorem 1 for operators whose kernel vanishes close to $\Sigma$. So far, working locally near a point which is not in $\Sigma$, we have only described the distribution of $h_j\partial_n e_i$ away from the singular set $\Sigma$. To get a complete description, we have to check that the $L^2$ norm can not accumulate at $\Sigma$. For this we are going to refine (3.1). Close to a point $x_0 \in \Sigma$ there exists $f_{i_1}$ and $f_{i_2}$ (possibly non unique) such that $f_{i_1}(x_0) = f_{i_2}(x_0) = 0$. Take $\chi \in C_0^\infty(\mathbb{R})$ equal to 1 close to 0 and consider (3.1) with $X$ replaced by

$$\chi\left(\frac{|f_{i_1}(x)| + |f_{i_2}(x)|}{\varepsilon}\right)X \tag{3.12}$$

and taking the limit $i \to +\infty$ (for fixed $\varepsilon > 0$) we obtain
(3.13)
$$\limsup_{k\to+\infty} \int_{\partial M} \chi\left(\frac{(f_{i_1}(x))^2 + (f_{i_2}(x))^2}{\varepsilon}\right)\left|h_k\frac{\partial e_k}{\partial n}\right|^2 \leq \langle \mu, \{|\xi|^2_{g(x)}, \chi\left(\frac{|f_{i_1}(x)| + |f_{i_2}(x)|}{\varepsilon}\right)X\}\rangle$$

But taking into account that $\mu$ is equal to $dL \otimes 1_{|\xi|=1}$, we see easily that the right hand side in (3.13) is bounded by $\mathcal{O}(\varepsilon)$ because the symbol is bounded by $1/\varepsilon$ but supported in a set of measure $\mathcal{O}(\varepsilon^2)$. As a consequence, passing to the limit $\varepsilon \to 0$, we obtain that the normal derivative can not accumulate at the singular set and consequently (3.11) holds for general semi-classical pseudodifferential operators on $\partial M$, tangential to $\Sigma$.

## 4. Robin and Neumann conditions

We consider in this section the eigenfunctions associated to the Robin or Neumann boundary conditions.

4.1. **A priori estimates.** By standard elliptic results (see for example [3, Appendix] in a slightly different context), we know that if $A$ is a semi-classical pseudo-differential operator supported in the elliptic region of the boundary, $\mathcal{E}$, one has for any $s \in \mathbb{R}$,

$$\|Ae_k^b\|_{H^s(\partial M)} \to 0 (k \to +\infty) \tag{4.1}$$

To deal with the hyperbolic region, we come back to (3.4). We obtain using the boundary condition

$$\int_M \frac{1}{ih_k}[h_k^2 \Delta_g - 1, \operatorname{Op}(a)(x, h_k D_x)]e_j \overline{e_j}|dx| \tag{4.2}$$
$$= -\int_{\partial M} a_2(x_n = 0, x', h_j D_{x'})(1 - R + h^2 K^2)(x', hD_{x'})e_j \overline{e_j}|dx'|$$

Remark that in (4.2) due to the elliptic regularity in the boundary (4.1) we can replace $K^2$ by $\Psi(x', h_k D_{x'}) K^2 \Psi(x', h_k D_{x'})^*$ with $\Psi$ equal to 1 for $|\xi'|_{g(x)} \leq 2$ and $\Psi$ equal to 0 for $|\xi'|_{g(x)} \geq 3$. The second order pseudodifferential operator $h^2 \Psi(x', D_{x'}) K^2 \Psi(x', D_{x'})^*$ can be viewed as a semiclassical operator of order 0 and principal symbol $|\Psi(x', \xi')|^2 \sigma_1(K)^2(x', \xi')$.

The left hand side in (4.2) is a second order differential operator in the $x_n$ variable. Using the relation (2.11) we see easily that in the left hand side of (4.2) we can eliminate the second order term in the normal derivative $h_k D_{x_n}$. Consequently the left thand side in (4.2) is bounded by $C_A(\|e_k\|_{L^2(M)} + h_k \|\nabla e_k\|_{L^2(M)}) \leq 2C(A)$. On the other hand, for any $q$ supported in the hyperbolic region, choosing $a_2$ such that

$$-a_2(x_n = 0, x', \xi')(1 - R + \sigma_1(K)^2)(x, \xi') = |q|^2 \geq 0 \tag{4.3}$$

we obtain, using sharp Gårding inequality (and the fact that by the trace Theorem we know that $\|e_k\|_{H^{-1/2}_{\text{loc}}(\partial M)} \leq C$ which gives a control on the remainder terms),

$$\|Qe_k\|_{L^2(M)} \leq C - \int_{\partial M} a_2(x_n = 0, x', h_j D_{x'})(1 - R + h^2 K^2)(x', hD_{x'})e_j \mid_{x_n = 0} \overline{e_j}|dx'| \leq C' \tag{4.4}$$

**Remark 4.1.** *Remark that the estimate above holds for any $k$, not only for a density-1 subsequence*

4.2. **The relationship between the interior and the boundary.** We come back to (4.2). As in section 3, the left hand side in (4.2) tends (for a density-1 subsequence) to

$$-\langle \mu, H_p(a) \rangle \tag{4.5}$$

from which we deduce, taking the limit in the right hand side

$$\lim_{j \in S \to +\infty} \int_{\partial M} a_2(x_n = 0, x', h_j D_{x'}, h_j)(1 - R + h^2 K^2)(x', hD_{x'})e_j \mid_{x_n = 0} \overline{e_j}|dx'| \tag{4.6}$$
$$= \langle \sqrt{1 - R} 1_{R \leq 1} |dx'||d\xi'|, a \times ((1 - R) + |\sigma_1(K)|^2) \rangle$$

which is Theorem 1 for the restricted class of symbols which can be written under the form

$$b(x', \xi') = a \times ((1 - R) + |\sigma_1(K)|^2) \tag{4.7}$$

(which is true if $b$ is equal to 0 close to $\mathcal{G}$).

## 5. The singular and glancing regions

So far we have obtained Theorem 1 for symbols which vanish close to the singular and glancing regions. To deal with these regions, we are going to use the following result obtained from heat kernel considerations in [8, Appendix] (not surprisingly since in this reference, the result was used for the same purpose):

**Proposition 5.1.** *For any $\varepsilon > 0$ there exists a function $\varphi_\varepsilon$ equal to 1 on $\Sigma$ but with support sufficiently close to $\Sigma$ such that*

$$(5.1) \qquad \left| \lim \frac{1}{N} \sum_{j=1}^{N} (\varphi_\varepsilon(x) e_j \mid_{\partial M}, e_j \mid_{\partial M})_{L^2(\partial M)} \right| = |\langle \varphi_\varepsilon \mid_{\partial M}, d\mu_b \rangle| \leq \varepsilon$$

**Proposition 5.2.** *For any $x_0 \in \partial M \setminus \Sigma$ and any $\varepsilon > 0$ there exists a symbol $\Psi_\varepsilon(x', \xi')$ equal to 1 on $\mathcal{G}$ close to $x_0$ but with support sufficiently close to $\mathcal{G}$ such that*

$$(5.2) \qquad \left| \lim \frac{1}{N} \sum_{j=1}^{N} \left( \Psi_\varepsilon(x', hD_{x'}) \Psi_\varepsilon(x', hD_{x'})^* e_j \mid_{\partial M}, e_j \mid_{\partial M} \right)_{L^2(\partial M)} \right| = |\langle d\mu_b, |\Psi_\varepsilon(x', \xi')|^2 \rangle| \leq \varepsilon$$

Indeed the equality in Proposition 5.1 is [8, Lemma 7.1], while the inequality follows if the support of $\varphi_\varepsilon$ is taken close enough to $\Sigma$; and Proposition 5.2 follows from [8, (11.14)] and the Karamata Tauberian Theorem, as in the proof of [8, Lemma 7.1].

To finish the proof of Theorem 1 for general operators, it is enough to prove that for any semi-classical pseudodifferential operator $A$ on $\partial M$ (tangential to $\Sigma$), one has

$$(5.3) \qquad \lim_{N \to +\infty} \frac{1}{N} \sum_{j=1}^{N} \left| \int_{\partial M} \mathrm{Op}(a)(x, \lambda_j^{-1} D_x, \lambda_j^{-1}) e_j^b \overline{e_j^b} d\sigma - \langle \mu_b, a(x, \xi) \rangle \right| = 0$$

By an easy argument, one deduces then *for fixed $A$* the existence of a density-1 subsequence such that (1.2) holds, then by a diagonal argument used in [10], [11] and [6] the existence of a density-1 subsequence such that (1.2) holds for any $A$ follows.

To prove (5.3), we use first Proposition 5.1 which allows to replace the operator $A$ by an operator $A_\varepsilon$ whose kernel is supported away from $\Sigma$, modulo an error of order $\mathcal{O}_A(\varepsilon)$. Consequently it is enough to prove the result for operators whose kernel are supported away from $\Sigma$; and from now on we suppose that $A$ satisfies this condition.

Write

$$\text{(5.4)} \quad \frac{1}{N}\sum_{j=1}^{N}\left|\int_{\partial M}\text{Op}(a)(x,\lambda_j^{-1}D_x,\lambda_j^{-1})e_j^b\overline{e_j^b}d\sigma - \langle\mu_b,a(x,\xi)\rangle\right|$$

$$\leq \frac{1}{N}\sum_{j=1}^{N}\Big|\int_{\partial M}\Big(\text{Op}(a)(x,\lambda_j^{-1}D_x) - \Psi_\varepsilon(x',\lambda_j^{-1}D_{x'})\text{Op}(a)(x,\lambda_j^{-1}D_x,\lambda_j^{-1})\Psi_\varepsilon(x',\lambda_j^{-1}D_{x'})^*\Big)$$

$$e_j^b\overline{e_j^b}d\sigma - \langle\mu_b,a(x,\xi)(1-|\Psi_\varepsilon(x',\xi')|^2)\rangle\Big|$$

$$+ \frac{1}{N}\sum_{j=1}^{N}\left|\int_{\partial M}\Psi_\varepsilon(x',\lambda_j^{-1}D_{x'})\text{Op}(a)(x,\lambda_j^{-1}D_x,\lambda_j^{-1})\Psi_\varepsilon(x',\lambda_j^{-1}D_{x'})^* e_j^b\overline{e_j^b}d\sigma\right|$$

$$+ \left|\langle\mu_b,a(x,\xi)|\Psi_\varepsilon(x',\xi')|^2\rangle\right|$$

Since the operator

$$\text{(5.5)} \quad \Big(\text{Op}(a)(x,\lambda_j^{-1}D_x,\lambda_j^{-1}) - \Psi_\varepsilon(x',\lambda_j^{-1}D_{x'})\text{Op}(a)(x,\lambda_j^{-1}D_x,\lambda_j^{-1})\Psi_\varepsilon(x',\lambda_j^{-1}D_{x'})^*\Big)$$

has its kernel supported away from $\Sigma$ and its symbol equal to 0 close to the glancing set $\mathcal{G}$, we obtain using (4.1) and (4.4) that

$$\text{(5.6)} \quad \left|\int_{\partial M}\text{Op}(a)(x,\lambda_j^{-1}D_x)e_j^b\overline{e_j^b}d\sigma - \langle\mu_b,a(x,\xi)\rangle\right|$$

is uniformly bounded with respect to $j$ and tends to 0 for a density-1 subsequence. Consequently the limit of the first term in (5.4) is equal to 0. On the other hand the contributions of the other terms in (5.4) are, according to Proposition 5.2, of size $\mathcal{O}(\varepsilon)$. This concludes the proof of (5.3).

Université Paris Sud, Mathématiques, Bât 425, 91405 Orsay Cedex
*E-mail address*: `Nicolas.burq@math.u-psud.fr`
*URL*: `http://www.math.u-psud.fr/ burq`